\documentclass{paper} 

\usepackage{graphicx}
\usepackage{amsbsy}
\usepackage{amsopn}

\title{Regularization of Kriging interpolation on irregularly spaced data.}

\author{Daniele Peri \footnote{\it{e-mail}: d.peri@iac.cnr.it}}
\institution{Istituto per le Applicazioni del Calcolo "M. Picone", \\ Consiglio Nazionale delle Ricerche, \\ Via dei Taurini, 19 - 00185 Roma - Italy}

\date{{\bf Keywords:} Kriging \sep Digital Twin \sep Artificial Intelligence \sep Interpolation}

\begin{document}

\maketitle

\begin{abstract}
Interpolation models are critical for a wide range of applications, from numerical optimization to artificial intelligence. The reliability of the provided interpolated value is of utmost importance, and it is crucial to avoid the insurgence of spurious noise. Noise sources can be prevented using proper countermeasures when the training set is designed, but the data sparsity is inevitable in some cases. A typical example is represented by the application of an optimization algorithm:  the area where the minimum or maximum of the objective function is assumed to be present is where new data is abundantly added, but other areas of design variable space are significantly neglected. In these cases, a regularization of the interpolation model becomes absolutely crucial. In this paper we are presenting an approach for the regularization of an interpolator based on the control of its kernel function via the condition number of the self-correlation matrix.
\end{abstract}

\section{Introduction}

In these years, the attention to the use of the so-called {\em digital twins} has gained a great popularity, particularly in the field of artificial intelligence (AI). The possibility to apply numerical methods replacing a complex and expensive real system, testing a number of different possibilities and configurations before producing the real object, is a great help for the design team. In addition, the comparison of the observed performances of a system and its ideal values, obtained as the output of the {\em digital twin}, can be useful for scheduling maintenance activities, thus reducing management costs.

Generally speaking, in most of the cases "{\em digital twin}" is a different way to call an interpolation model, since the response of the {\em digital twin} is commonly based on the manipulation of a number of observed values of the different quantities characterizing the system status: data interpolation represents therefore the core element. If we need high quality in the data interpolation, and of the {\em digital twin} consequently, a great attention must be dedicated to the design of the training set, uniformly distributing the sample points. Unfortunately, if the building of the interpolator is not static, and new points can be added dynamically, the resulting training set could lost its uniformity, and also the quality of the interpolation could become poor.

An example is represented by the use of the {\em digital twin} inside an optimization loop: classically, the optimization algorithm is using the {\em digital twin} instead of the true objective function (if computationally expensive) for the detection of the optimal solution. Once the current optimal solution is detected, the candidate solution is verified and a new point can be added to the training set. If this operation is repeated, it is common that all the new points will belong to the same area in the design variable space, creating an imbalance in the point distribution.

The {\em digital twin} can be also based on real (non simulated) data, collected when the system is operational, or can be obtained from measurements of a natural phenomenon.
Irregularly-sampled time series occur in many domains including healthcare, geophysics and many more. They can be challenging to model because they do not naturally yield a fixed-dimensional representation as required by many standard machine learning model. Adding new points in these cases could enhance the difficulties, because new data are not completely under the control of the researchers, so they are often not improving t database under the standpoint of its uniformity and regularity.

To mitigate the effect of the lack of uniformity of the samples, a localization of the interpolation may be performed: the {\em digital twin} can divide the entire database into different parts, accessing the most appropriate one depending on the required configuration, but some irregularities at the borders of the splitting areas could appear. {\em Viceversa}, if we prefer to preserve the integrity of the database, producing a single interpolator for the complete set of possible system configurations, a kind of regularization technique is mandatory.

Some solutions have been already presented in literature, and their number is quite large. The proposed approaches can be categorized in two great families: the so called {\em dynamic kriging}, presented and reviewed in \cite{Zhao2011}, and more generic regularization methods (an example is presented in \cite{Zhang2019}).

The dynamic kriging tunes the optimal values of the kernel parameters by optimizing a number sub-problems resulting from the cross-validation technique: it has been already preliminarily tested in \cite{Peri2015}. The procedure could be extremely time consuming, since the number of subproblem to be solved could be very large, depending on the number of available sampling points, and the final parameters are obtained solving a suite of sub-problems all different than the final one. Anyway, this technique has been demonstrated to be pretty efficient.

In the regularization methods, the kernel function is modified to include the regularization terms: in this case, the solved problem is the same as the final one, but the use of regularization terms can transform kriging from interpolation to approximation, losing precision locally.

In this article we present a regularization method suitable for a generic kernel method, or in any case for all methods that require the inversion of a square matrix for the determination of the key elements of the interpolator (typically, the weights of some kernel functions). There is not a modification of the kernel function e.g. by using some penalties term, since we are seeking for the optimal value of the parameters ruling the kernel functions. The full training set is used, without any further manipulation or reduction, and all the available information are exploited. Results are presented for the kriging interpolation method, that is responding to the aforementioned requirements. Numerical evidences are absolutely satisfactory under the perspective of preserving a regular behavior of the interpolation, increasing the overall credibility of the produced output.

\section{Basics about kriging}

Since we are producing results using the kriging method, we are going here to recall some basic elements. Among the different techniques for data interpolation, kriging represents a very attractive option. It was firstly introduced in \cite{Matheron1963}, and since then a large number of studies and variants have been proposed. A rather complete review paper is \cite{Kleijnen2009}. Following the notation reported in \cite{Kleijnen2009}, the estimated value $y(\pmb{d})$ can be obtained as

\begin{equation}\label{Pesi}
y(\pmb{d}) = \pmb{\lambda}(\pmb{d},\pmb{D})^\prime \pmb{w}(\pmb{D})
\end{equation}

where {\bf d} indicates the position of the point where the estimate is required, {\bf D} is the position of the training points and {\bf w(D)} is the value of the function to be interpolated at the sample points. 

To select the optimal values for the weights $\pmb{\lambda}$, a least-squares problem is solved minimizing the Sum of Squared Residuals, that is, minimizing the Mean Squared Error (MSE) of the predictor y(d)

\[
 \operatorname*{min}_{\lambda}\mathsf{MSE}= \operatorname*{min}_{\lambda}\big[ E(y(\mathrm{d})-\mathsf{w}\mathsf{(d))}^{2} \big]
\]

where {\bf d} may be any point in the design variable space. Moreover, this minimization must account for the condition that the predictor is unbiased:

\[
E(y(\pmb{d})) = E(w(\pmb{d})) 
\]

where in deterministic simulation E(w({\bf d})) may be replaced by w({\bf d}). It can be proven that the solution of this constrained minimization problem implies that the sum of the weights is 1. Finally, it can be also proven that the optimal value of the weights is

\begin{equation}\label{KrigingSys}
\lambda_{o}=\pmb{\Gamma}^{-I}\left[\gamma+\pmb{1}\,\frac{1\,-\,\pmb{1}^{\prime}\pmb{\Gamma}^{-1}\gamma}{\pmb{1}^{\prime}\pmb{\Gamma}^{-1}\pmb{1}}\right]
\end{equation}

where $\pmb{\Gamma} = cov(w_i,w_{i^\prime})$ is the covariance matrix involving the training set data, and 
$\pmb{\gamma} = cov(w_i,w_0)$ is the covariance vector including the training set data and the interpolation location. See \cite{Kleijnen2009} for details.

What is still missing is the definition/selection of the covariance function. In the original formulation, the concept of {\em semi-variogram} is adopted instead, and it can be estimated experimentally using the training set\cite{Matheron1963}. Unfortunately, the {\em semi-variogram} represents a statistical quantity, and the number of training points is rarely adequate in real applications for producing a significative output, so that some analytical function are used instead. Practically, an algebraic expression, typically an exponential law, is used to fit the (few points of the) sample points, obtaining a kind of experimental semi-variogram: In fact, the data do not clearly indicate a trend, and even less they are able to indicate any continuous curve to apply: consequently, for any algebraic model adopted, the difference between the experimental values and the semi-Approximate variogram is usually large. A commonly adopted function is

\begin{equation}\label{Vario1}
\rho(\mathrm{h})=\exp\left[-\sum_{j=1}^{k}\theta_{j}h_{j}^{p_{j}}\right]=\prod_{j=1}^{k}\exp[-\theta_{j}h_{j}^{p_{j}}]
\end{equation}

$\theta_j$ denotes the importance of input $j$: the higher $\theta_j$ is, the less effect input $j$ has. $p_j$ denotes the smoothness of the correlation function: $p_j$ = 2 implies an infinitely differentiable function. As a consequence, the correlation function is a function of these two parameters, to be selected.

The tuning of the free parameters, $(\theta_{j},p_{j})$ if the previous formulation for the covariance function is adopted, can be performed splitting the available samples in two different sets: {\em training} and {\em validation}, using the first one to derive $\pmb{\lambda}$ and the second one to tune $(\theta_{j},p_{j})$. An example is reported in \cite{Peri2015}. But this technique is difficult to apply in practice, since the number of sampled points is generally not large, and the use of some of these points as {\em validation} points is further impoverish the training set. By the way, the exactness of the solution is guaranteed by kriging on the training point, so that the shift of a training point in the validation set is not intrinsically beneficial.

In the following, we are going to fix $p_j$ = 2, so that the only free parameters are $\theta_{j}$.

\section{Kriging regularization - numerical experiments}

The basic idea behind the present scheme for the regularization of kriging is simple. The determination of the weights $\lambda_j$ is obtained by solving a linear system: if we find a way to maximize the condition number of the matrix of the linear system, the solution could be much more stable, and possibly much more regular. To verify this hypothesis, using the aforementioned definitions, an optimization problem has been solved, minimizing the condition number of the self-correlation matrix $\pmb{\Gamma}$ acting on the parameters $\theta_{j}$. The starting point for the optimization problem is obtained by investigating a number of alternative solutions for $\theta_{j}$, producing a perturbation around an initial value. Starting from the so determined best solution (best condition number), a local minimization algorithm of the class of the Direct Search Methods \cite{Powell2006} is applied for the fine tuning of the parameters $\theta_{j}$.

Six different 2D test functions have been adopted. 2D functions have been selected in order to provide a clearer evidence of the obtained results: the method is obviously valid for an arbitrary space dimension. The selected six functions are:

\begin{enumerate}
\item Griewank function \\ $f(x)=\sum_{i=1}^{d}{\frac{x_{i}^{2}}{4000}}-\prod_{i=1}^{d}\cos\left({\frac{x_{i}}{\sqrt{i}}}\right)+1$ \\
 where $x_i \in [-5:5]$
\item Sasena function \\ $f(x)=2 + 0.01 (x_2-x_1^2)^2 + (1 -x_1)^2 + 2 (2 -x_2)^2 + 7 sin(0.5 x_1) sin(0.7 x_1 x_2)$ \\
 where $x_i \in [0:5]$
\item Franke's function \\ $f(x)=0.75\exp\left(-\frac{(9x_1-2)^2}{4}-\frac{(9x_2-2)^2}{4}\right)+0.75\exp\left(-\frac{(9x_1+1)^2}{49}-\frac{9x_2+1}{10}\right) + 0.5\exp\left(-\frac{(9x_1-7)^2}{4}-\frac{(9x_2-3)^2}{4}\right)-0.2\exp\left(-(9x_1-4)^2-(9x_2-7)^2 \right)$ \\
 where $x_i \in [0:1]$
\item G-function \\ $f(x)=\prod_{i=1}^{d}\frac{\left|4x_{i}-2\right|+a_{i}}{1+a_{i}},~{\mathrm{where}} \, a_{i}={\frac{i-2}{2}},~\mathrm{for~all~}i=1,\cdots,d$ \\
 where $x_i \in [0:1]$
\item Irregular TF \\ $f(x)=exp(x_1)/5-x_2/5+x_2^6/3+4 x_2^4-4 x_2^2+7/10 x_1^2+x_1^4+3/(4 x_1^2+4 x_2^2+1)$ \\
 where $x_i \in [-1:1]$
\item Cosin2 \\ $f(x)=cos(10 x_1)+sin(10 x_2) + x_1 x_2$ \\
 where $x_i \in [0:1]$
\end{enumerate}

\begin{figure}[htb]
\centering
\includegraphics[width=\textwidth]{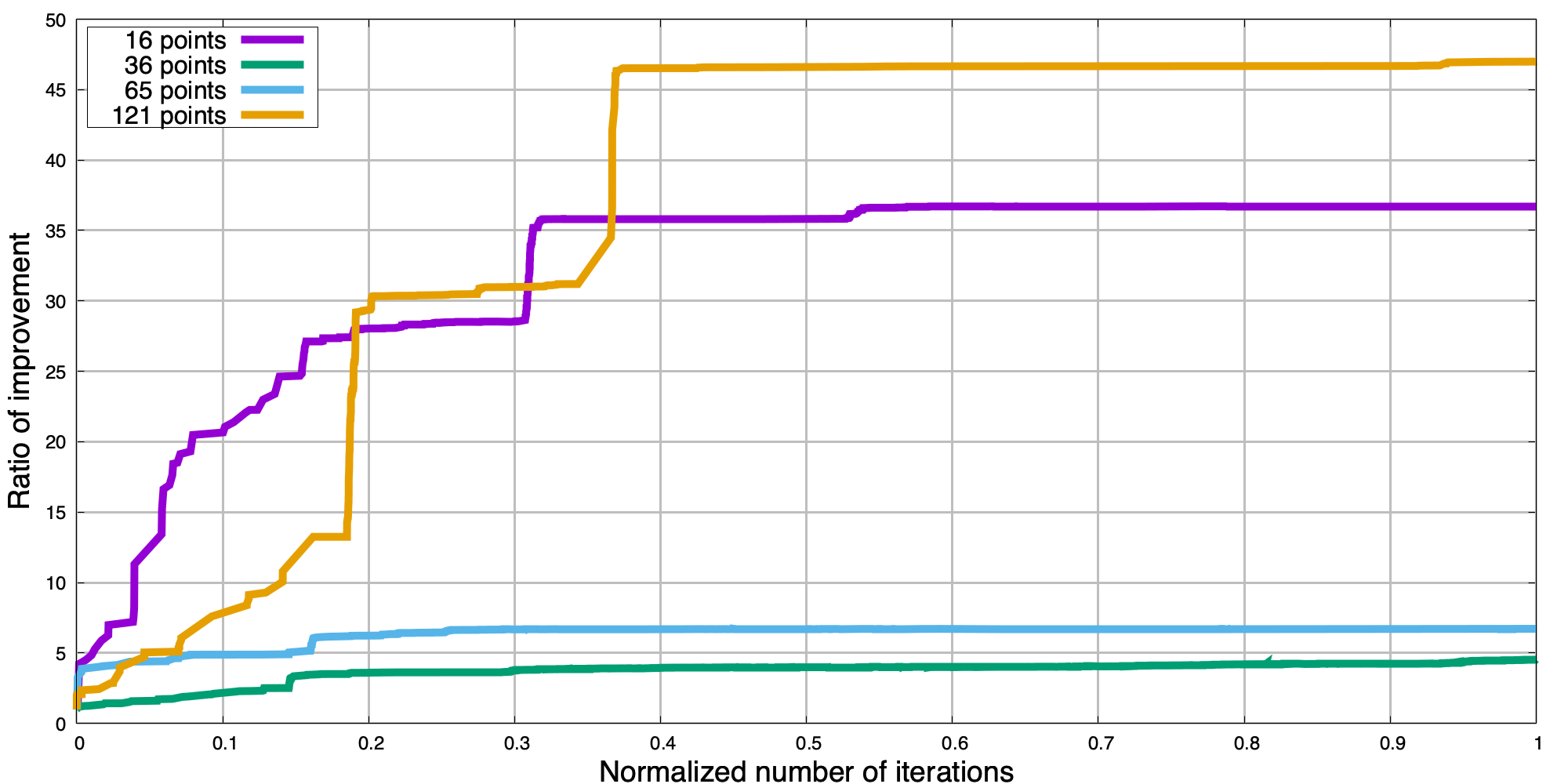}
\caption{Convergence of the condition number of the self-correlation matrix for a variation of the semi-variogram coefficients $\theta_j$. The number of iterations and condition number are normalized respectively by the maximum number of iterations and the initial value of the condition number. 4 different density of the training points, randomly generated. {\em Franke's} function.
        }
\label{fig:Conve}
\end{figure}

A first comment comes from the observation of figure \ref{fig:Conve}, where four different examples, with different number of training points, are compared. The improvement of the condition number of the self-correlation matrix, defined as the ratio between the current value and the initial value, is reported as a function of the non-dimensional value of the iterations of the optimization problem (the ratio between the current number of iteration and the maximum value). It is clear how, for a moderate number of training points, the improvement of the condition number is not large, while it becomes high in case of a very small or very large number of training points. The reasons for this improvement are exactly opposite: while for 16 training points the absolute value of the condition number is already high, a great improvement achieved in relative terms translates into a small change of the model in practice. On the opposite, the condition number can be largely improved in absolute value when the training set is large and unevenly distributed, and we have both a great improvement in absolute and relative terms. Figure \ref{fig:Conve} is demonstrating that very good results that have been obtained in terms of condition number.

Due to the nature of the kriging method, the matrix $\pmb{\Gamma}$ is depending only from the correlation function, and not from the shape of the objective function. As a consequence, if the point distribution, in a normalized design variable space, is the same, also the results for different objective functions are exactly the same. For this reason, the figure \ref{fig:Conve} is here reported for the case of the {\em Franke's} function, but in reality it is valid also for any other objective function. The results for all the algebraic functions listed above are shown below, but the use of one function instead of another is completely equivalent, unless the training set distribution is not changed. 

\begin{figure}[htb]
\centering
\includegraphics[width=0.99\textwidth]{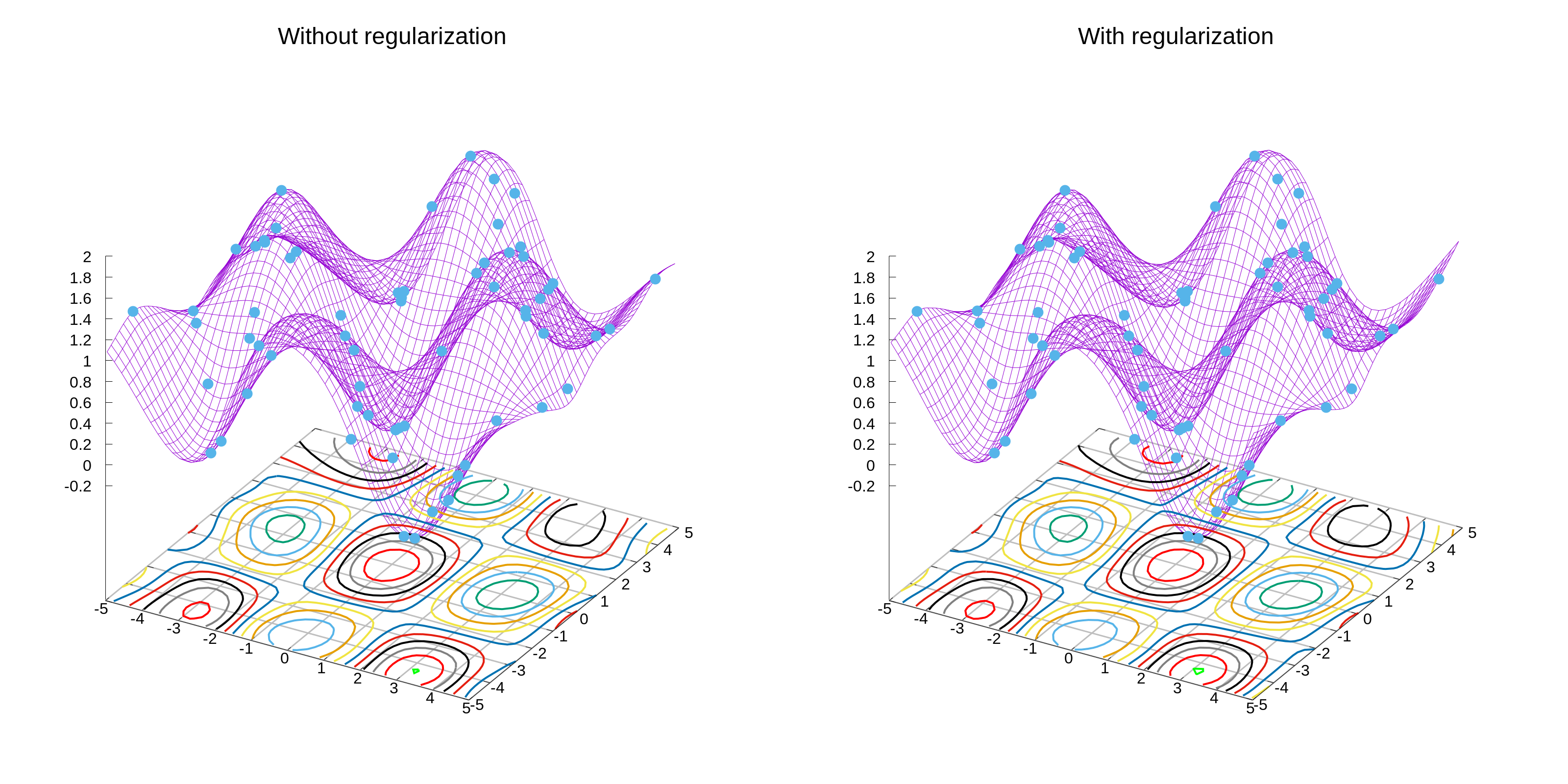} \\
\includegraphics[width=0.99\textwidth]{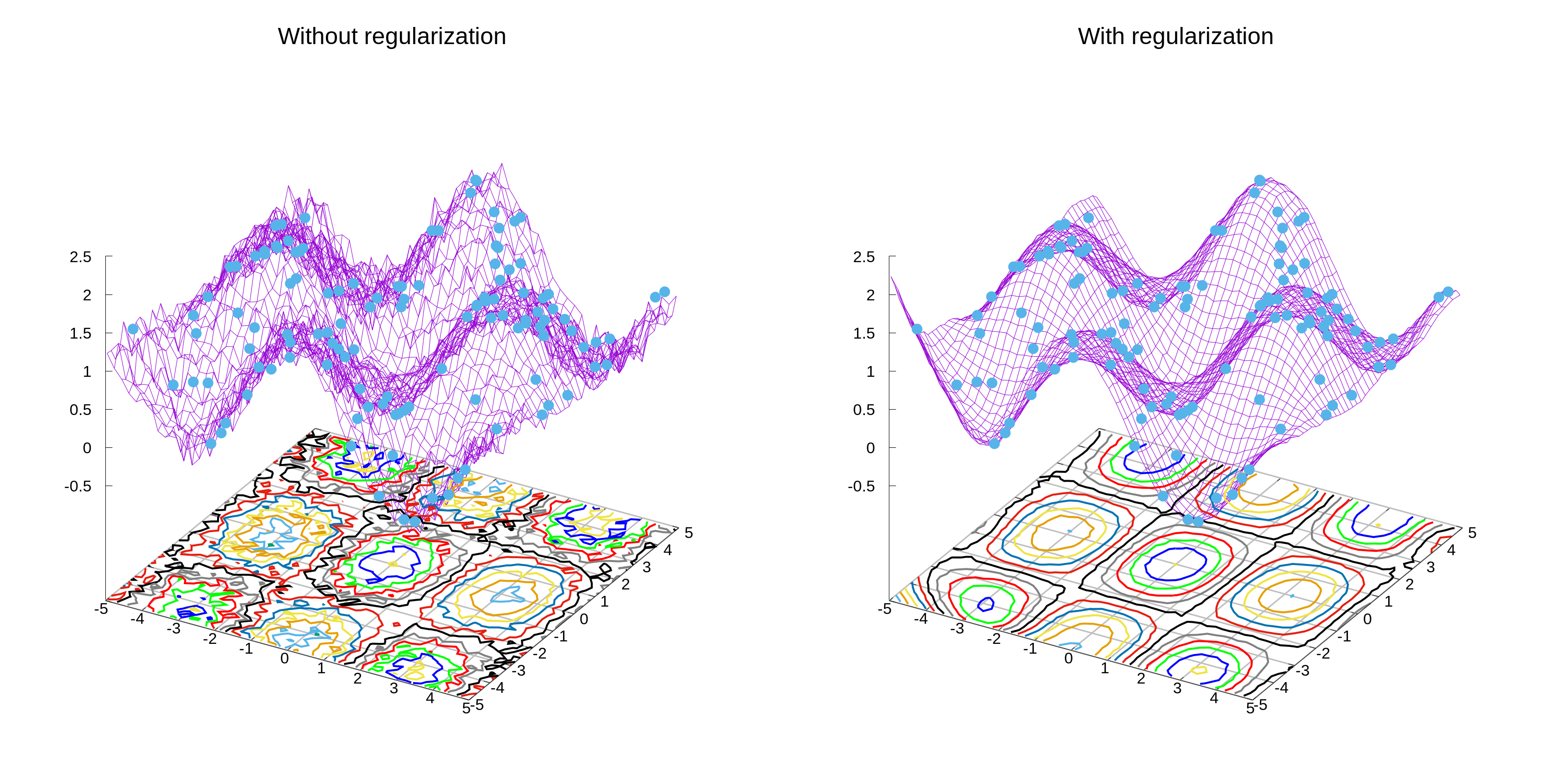}
\caption{Comparison between the interpolation obtained for an increasing number of irregularly-generated (random) points, with and without regularization. Two-dimensional Griewank function. 64 (on top) and 121 (on bottom) training points. Left without regularization, right with regularization.
        }
\label{fig:Griewank}
\end{figure}

Results about the regularization of the interpolating function are reported in figure \ref{fig:Griewank}: only the graphical results for last two training sets are reported, since for smaller values of the training points differences are not visible. On top of figure \ref{fig:Griewank}, the case with 64 training points is reported, and there is not a clear difference between the case with and without regularization. On the contrary, for the case with 121 training points, reported in the bottom part of figure \ref{fig:Griewank}, the difference is clear. To highlight irregularities, contour levels are provided at the bottom of the graph. In the case without regularization, the presence of high-frequency noise is clear, while for the regularized interpolation we have a global behavior comparable with the case with 64 training points, but with a sharper distinction between the basins of attraction. Summarizing, we can affirm that the regularization procedure produces beneficial effects.

\begin{figure}[htb]
\centering
\includegraphics[width=0.49\textwidth]{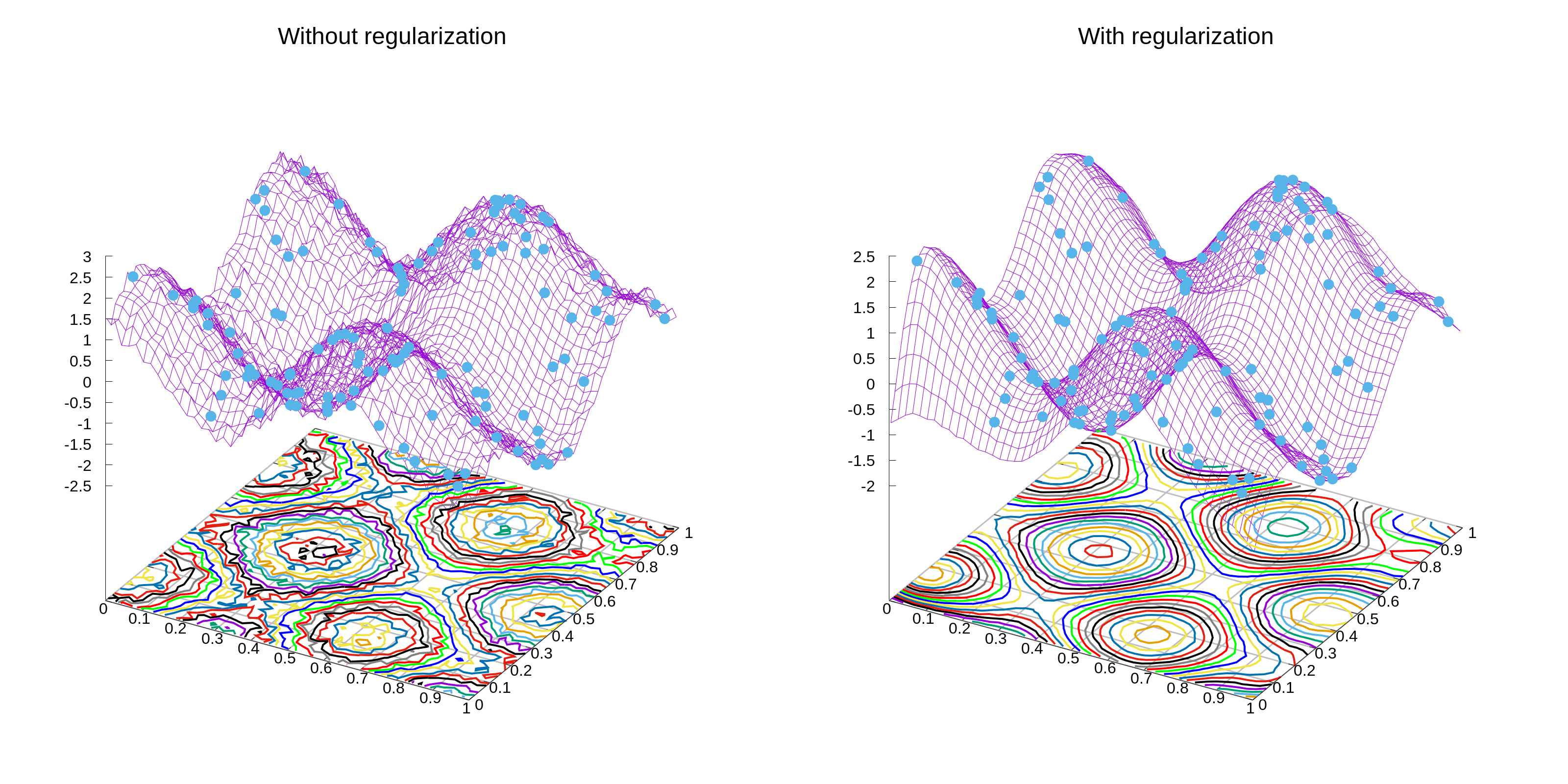}
\includegraphics[width=0.49\textwidth]{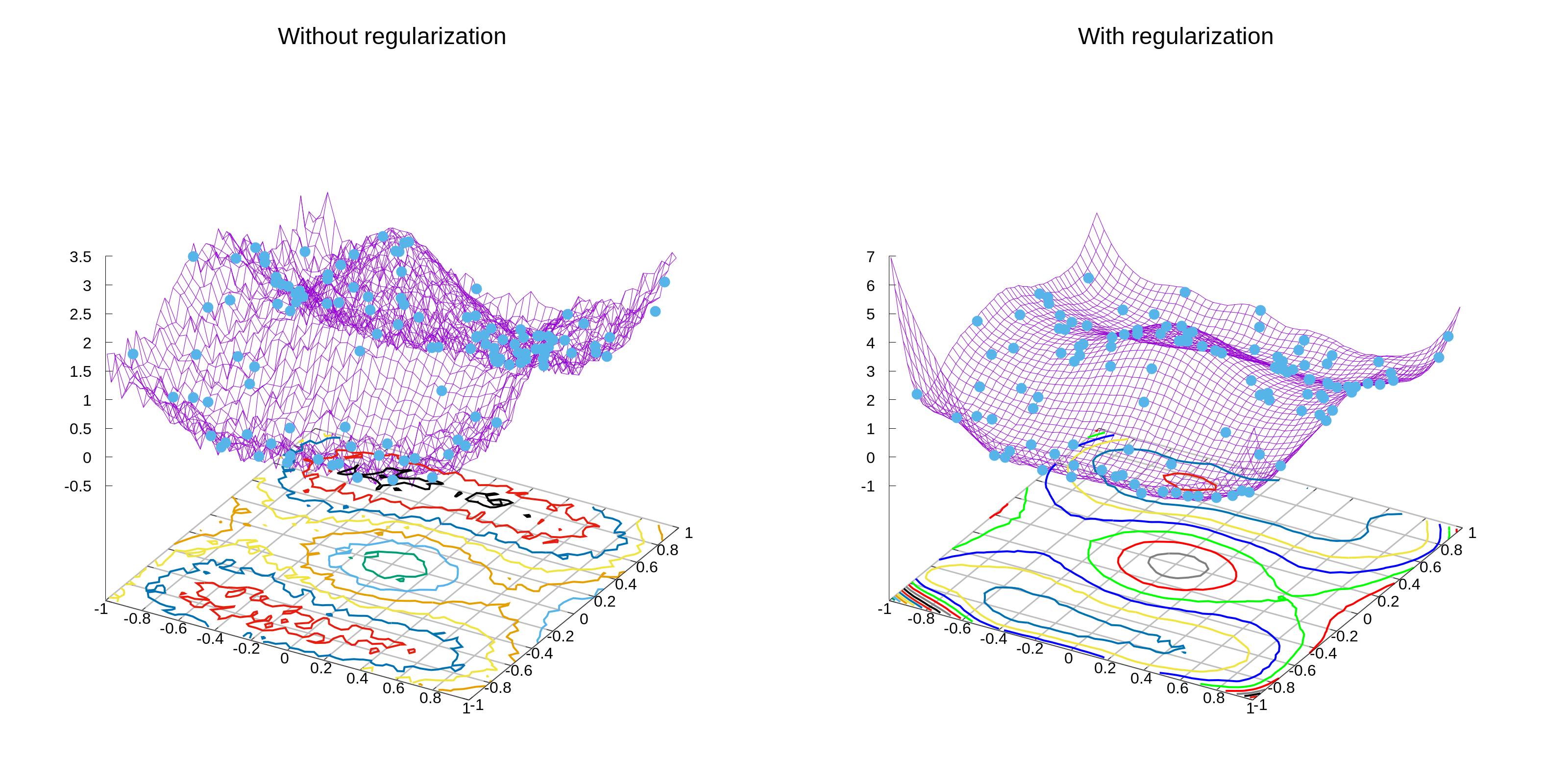} \\
\includegraphics[width=0.49\textwidth]{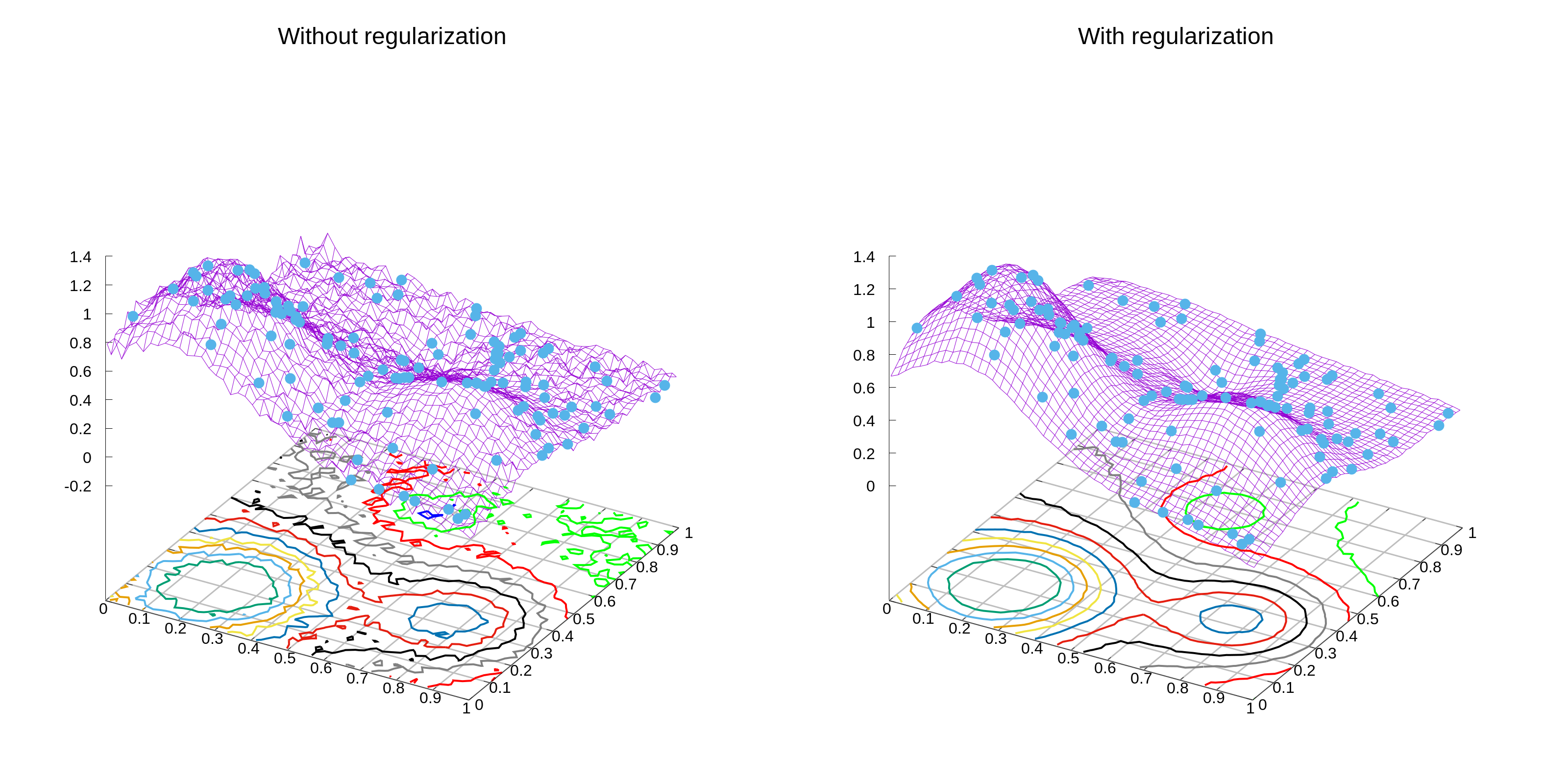}
\includegraphics[width=0.49\textwidth]{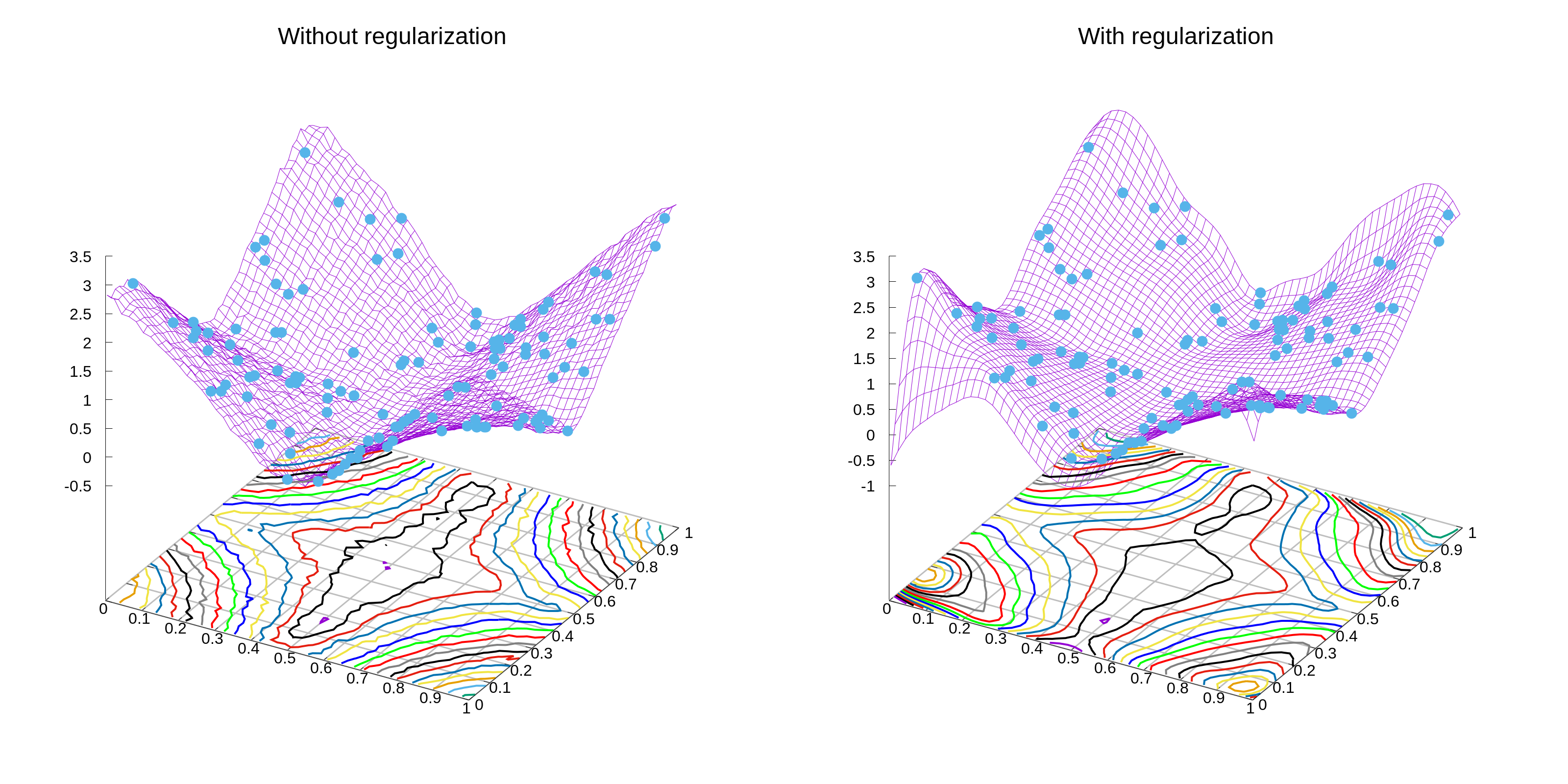} \\
\includegraphics[width=0.49\textwidth]{gri121.png}
\includegraphics[width=0.49\textwidth]{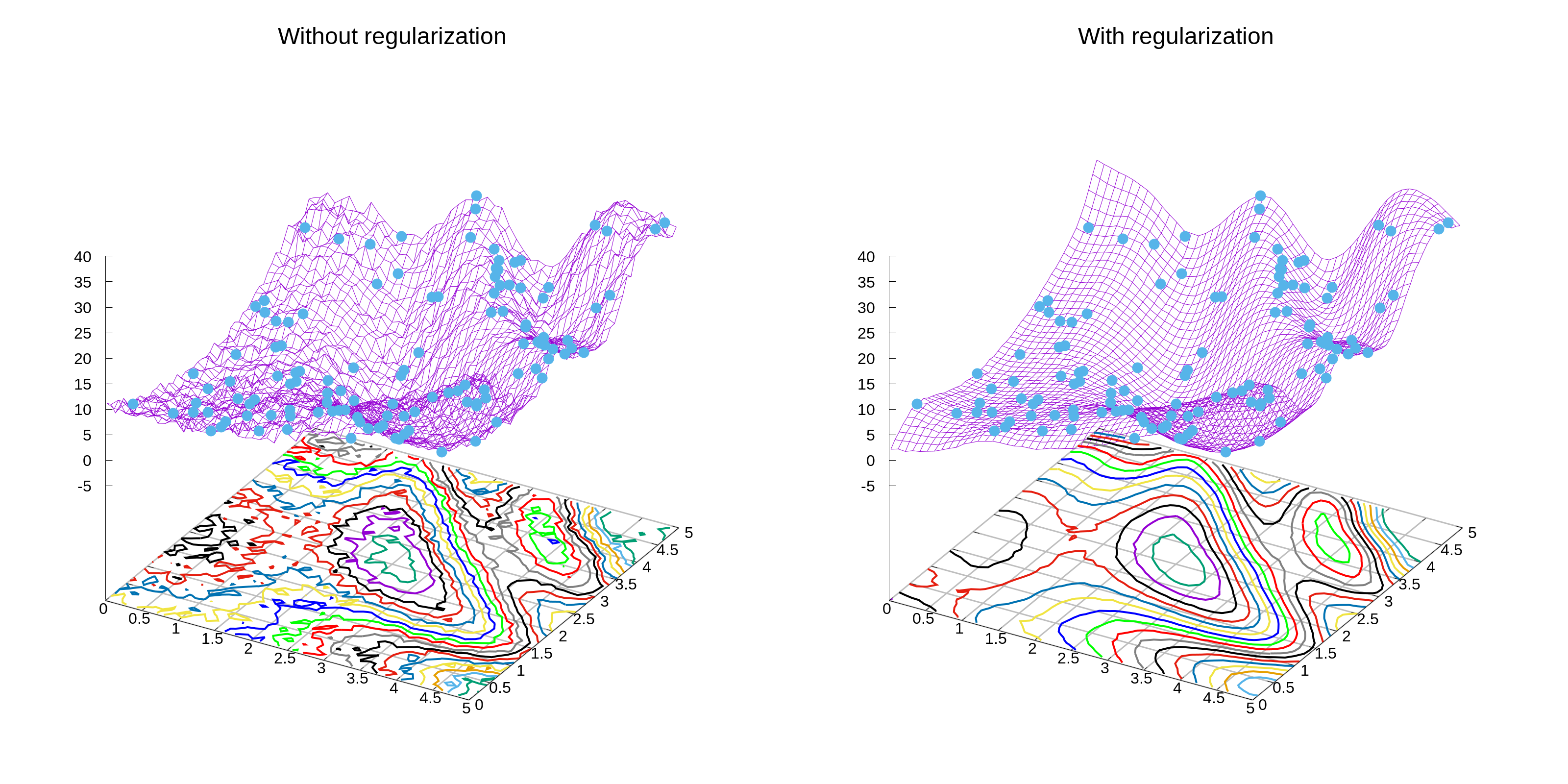}
\caption{Comparison between the interpolation obtained for 121 randomly selected training points with and without regularization. Two-dimensional test functions. From left to right, top to bottom: Cosin2, Irregular TF, Franke's function, G-function, Griewank function, Sasena function. For each function, the two global reconstructions (without and with the regularization algorithm) are reported.
        }
\label{fig:All}
\end{figure}

The same behavior is observed, as anticipated earlier, for all the six investigated algebraic function reported in figure \ref{fig:All}. The noise is clearly eliminated by the regularization procedure, improving the overall quality and consistency of the interpolation: this feature is of paramount importance for the use of kriging in combination with local optimization algorithms.

\begin{figure}[htb]
	\centering
	\includegraphics[width=0.99\textwidth]{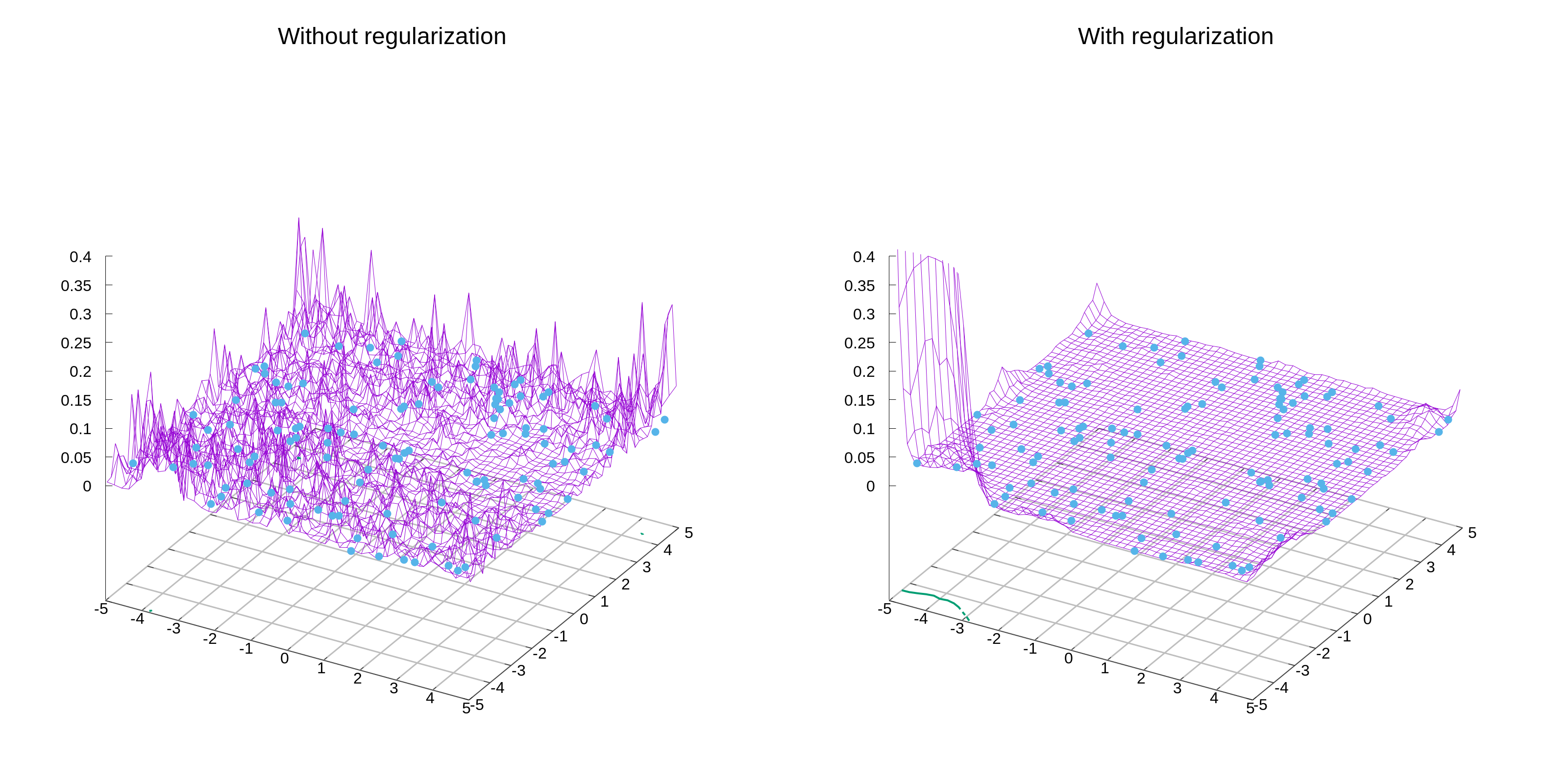} \\
	\includegraphics[width=0.99\textwidth]{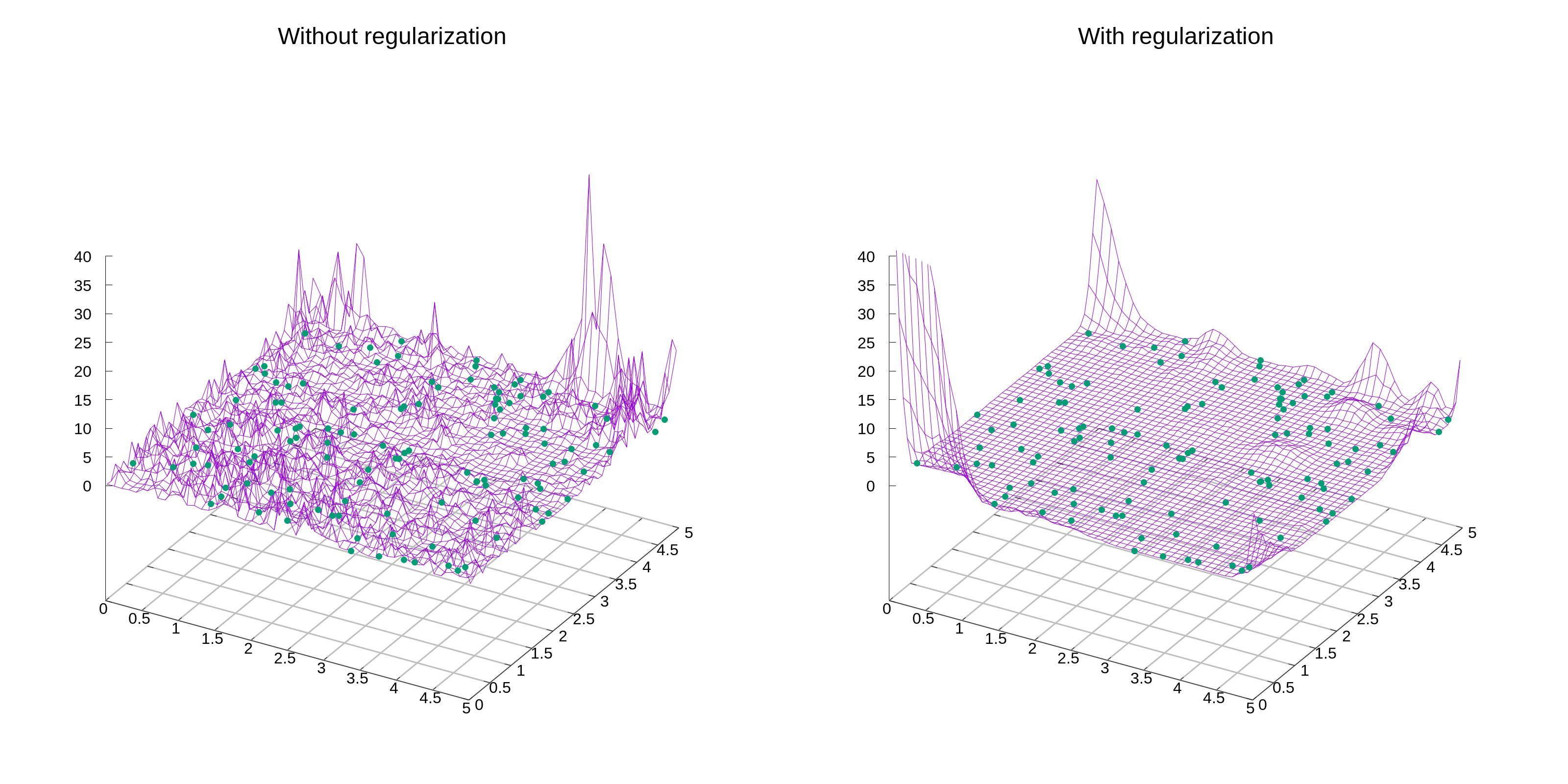}
	\caption{Comparison between the interpolation errors for the Griewank function (on top) and the Sasena function (on bottom) without (on left) and with (on right) the use of the proposed regularization algorithm.
	}
	\label{fig:GriewankErr}
\end{figure}

Figure \ref{fig:GriewankErr} provides a local quantification of the differences between the interpolator and the true value of the function for the {\em Griewank} and {\em Sasena} functions. Two examples are given, even if not necessary, simply to allow to notice the great similarities: minor differences are linked only to the different absolute value of the functions, and the location of the areas where the discrepancies are greatest is absolutely the same. The dots at the base of the graph are reporting the position of the training points. We can see how regularization reduces the differences to almost zero in the largest part of the domain. On the contrary, in a corner, where training is poor, regularization produces a local increase of the  differences: it is quite obvious that regularization cannot improve the content of information, but only reshape interpolation more regularly, so that regions where the point density is low are subject to these consequences.
        
\section{Conclusions}

A regularization method for a generic kernel interpolator, based on the observation of the condition number of the kernel matrix, has been here described and successfully tested. Implementation and results has been produced for the kriging interpolator only, but we can use this algorithm also for different kernel methods: this is desirable for a larger comparison, in the future.
 
Also a direct and systematic comparison with other regularization methods, e.g. the ones based on the split of the sampling points in two different sets, {\em training} and {\em validation}, as reported in \cite{Zhang2019} and \cite{Peri2015}, could be also desirable.

\section*{Acknowledgements}

D.P. research was partly funded by Italian Minister of Instruction, University and Research
(MIUR) to support this research with funds coming from PRIN Project 2017 (No. 2017KKJP4X
entitled "Innovative numerical methods for evolutionary partial differential equations and
applications").

\section*{Conflict of interest}

The author declares no conflicts of interest.

\section*{Data avaliability}

The data that support the findings of this study are available on request from the author.



\begin{thebibliography}{99}

\bibitem[{Kleijnen(2009)}]{Kleijnen2009}
Kleijnen, J.~P., 2009: Kriging metamodeling in simulation: A review.
  \textit{European Journal of Operational Research}, \textbf{192~(3)},
  707--716.

\bibitem[{Matheron(1963)}]{Matheron1963}
Matheron, G., 1963: Principles of geostatistics. \textit{Economic Geology},
  \textbf{58~(8)}, 1246--1266.

\bibitem[{Peri(2015)}]{Peri2015}
Peri, D., 2015: Improving predictive quality of kriging metamodel by variogram
  adaptation. \textit{MARINE VI : proceedings of the VI International
  Conference on Computational Methods in Marine Engineering}, CIMNE, Ed.,
  Vol.~1.

\bibitem[{Powell(2006)}]{Powell2006}
Powell, M. J.~D., 2006: \textit{The NEWUOA software for unconstrained
  optimization without derivatives}, 255--297. Springer US, Boston, MA.

\bibitem[{Zhang et~al.(2019)}]{Zhang2019}
Zhang, Y., W.~Yai, and X.~Chen, 2019: A regularization method for constructing
  trend function in kriging model. \textit{Structural and Multidisciplinary
  Optimization}, \textbf{59}, 1221--1239.

\bibitem[{Zhao et~al.(2011)}]{Zhao2011}
Zhao, L., K.~K. Choi, and I.~Lee, 2011: Metamodeling method using dynamic
  kriging for design optimization. \textit{AIAA Journal}, \textbf{49},
  2034--2046.

\end{thebibliography}
\end{document}